\documentclass[12pt, a4paper]{amsart}

\usepackage[all]{xy}
\usepackage{amsmath}
\usepackage{amssymb}
\usepackage{amsfonts}
\usepackage{graphicx}
\usepackage[mathscr]{eucal}
\usepackage{enumerate}
\usepackage[latin1]{inputenc}

\theoremstyle{definition}
\newtheorem{definition}{Definition}[section]
\theoremstyle{remark}

\newtheorem*{remark}{\textbf{Remark}}

\newtheoremstyle{estilo}{\topsep}{\topsep}{\slshape}{}{\bfseries}{.}{ }{}
\theoremstyle{estilo}
\newtheorem{theorem}[definition]{Theorem}
\newtheorem{corollary}[definition]{Corollary}

\newtheorem{proposition}[definition]{Proposition}
\def\ot{\otimes}
\def\op{\oplus}

\renewenvironment{proof}{{\noindent\sc Proof\;}}{\qed\\}

\DeclareMathOperator{\Id}{Id} \DeclareMathOperator{\End}{End}
\DeclareMathOperator{\Hom}{Hom} \DeclareMathOperator{\im}{Im}

\newcommand{\dtext}[1]{\emph{\textbf{#1}}}

\def\lto{\longrightarrow}

\DeclareMathOperator{\Dim}{dim} \DeclareMathOperator{\Char}{char}

\title{On the classification and properties of noncommutative
duplicates}
\author{Javier L\'{o}pez Pe\~{n}a}
\author{Gabriel Navarro}
\address{Department of Algebra,
University of Granada\\
Avda. Fuentenueva s/n, E-18071, Granada, Spain}
\email{jlopez@ugr.es, gnavarro@ugr.es}
\thanks{This research has been partially supported
by the projects MTM2004-08125, MTM2004-01406 and FQM-266 (Junta de Andaluc\'{i}a
Research Group). First author is partially supported by Spanish MEC-FPU
grant AP2003-4340.}

\date{}

\begin{document}
\maketitle

\begin{abstract}
    We give an explicit description of the set of all factorization structures,
or twisting maps, existing between the algebras $k^2$ and $k^2$, and
classify the resulting algebras up to isomorphism. In the process we
relate several different approaches formerly taken to deal with this
problem,
filling a gap that appeared in a recent paper
\cite{cibils}. We also provide a counterexample to a result
concerning the Hochschild (co)homology appeared in \cite{guccione}.

\end{abstract}

\section*{Introduction}
Interest on studying the structure of twisted tensor products of
algebras comes from different sources. On the one hand, they are in
one to one correspondence with factorization structures, defined in
\cite{majid}, on the other hand, they admit a reinterpretation as a
representative for a cartesian product of spaces in the framework of
noncommutative geometry (cf \cite{cap}, \cite{jara}). Though most of
the purely algebraic properties of this construction have a very
strange behaviour, which has meant a serious pitfall for most trials
of systematically describing properties of twisted tensor products,
the geometrical interpretation of this structure provides quite a
good intuition, becoming the approach which has proved to be the
most successful in the recent years.

Whilst there certainly exists a strong motivation in different areas
of algebra for studying the structure of twisted tensor products,
there is one basic problem concerning them which turns out to be of
much more fundamental nature. Namely, the classification of all
different twisted tensor products that we can obtain starting from
two given algebras, say $A$ and $B$. If we are dealing with unital
algebras, the equivalence between the existence of a twisted tensor
product structure and the existence of a twisting map simplifies the
problem into fixing a pair of vector spaces, $A$ and $B$, and
finding out linear maps between $B \otimes A$ and $A \otimes B$
satisfying certain properties. In some particular situations, for
instance when both $A$ and $B$ are finite dimensional, we can easily
rewrite the twisting conditions in terms of the matrix elements of
the linear map (once we have fixed some bases in the vector spaces),
thus obtaining some polynomial equations that the component of a
matrix must satisfy in order to be the representation of a twisting
map in the given bases. As a consequence, we may look at the set of
twisting maps between $A$ and $B$ as an affine subvariety of
$M_{m\times n}(k)$, where $m$ is the dimension of $A$ and $n$ is the
dimension of $B$. The classification of the twisting maps is
therefore equivalent to the description of this algebraic variety.
Some steps on describing the variety of twisting maps has been given
by Cibils in \cite{cibils}.

The second obvious problem that we face when trying to classify the
twisted tensor products between two algebras is the fact that
different twisting maps may give rise to isomorphic algebra
structure. In terms of the variety of twisting maps, this problem
boils down to study a quotient variety of our original one.
Unfortunately, so far no groundwork that can simplify the
isomorphism problem is known, and we are bond to deal with each case
separately.

Though several attempts to describe the structure and behavior of
twisted tensor products of algebras has been taken, the inherent
difficulty of the task has caused a number of flaws in different
papers on the topic. Our purpose in this work is to clarify, through
the exhaustive description of the space of twisting maps between two
particular algebras and the orbitspace of isomorphism classes
associated to them, some points that have recently appeared in the
literature under an incorrect form. In particular, we fix a subtle
mistake in Cibils' classification of noncommutative duplicates. We
also give a counterexample to several results given by J.A. Guccione
and J.J. Guccione in \cite{guccione}, concerning Hochschild homology
of twisted tensor products.

It is worth noting that our example stress the fact that the
homological properties of a twisted tensor product is usually not
well behaved with respect to the same properties in the factors,
even if we impose some conditions on the twisting map. However, if
we impose some stronger conditions on the involved algebras, it is
still possible to recover some regular behaviour. In particular,
Koszulity of the algebras is enough to guarantee Koszulity of the
twisted tensor product, giving rise to a K\"{u}nneth formula
describing the homology of the product algebra in terms of the
homology of the factors. Further details in this matter can be found
on \cite{jara2}.

\section{Preliminaries}\label{preliminaries}

Along the present work, $k$ will stand for a ground field, which we
will require in some situations, explicitly mentioned in the text,
to be of characteristic different from 2. All algebras will be
supposed to be associative, unital, $k$--algebras. The term
\emph{linear} will always mean \emph{$k$--linear}, and the unadorned
tensor product $\otimes$ will stand for the usual tensor product
over $k$. For an algebra $A$ we will write $\mu_A$ to denote the
product in $A$ and $u_A:k\rightarrow A$ its unit. We will also
identify every object with the identity map defined on it, so that
$A\otimes f$ will mean $\Id_A\otimes f$.

The notion of twisted tensor product has been independently discovered a number
of times. Some of the earliest papers giving the definition that we shall use
are \cite{cap}, \cite{cartier}. Let $A$, $B$ be two algebras, and $\tau:B\ot
A\to A\ot B$ a linear map. We say that $\tau$ is a \dtext{twisting map} if it
satisfies the following conditions:
\begin{eqnarray}
    & \tau(b\ot 1)=1\ot b,\quad \tau(1\ot a)=a\ot 1, \quad \forall \; a\in A,
\;b\in B, \label{tw1}\\
    & \tau \circ (B\ot \mu _A)=(\mu _A\ot B) \circ (A\ot \tau)\circ (\tau\ot
A), \label{tw2}\\
    & \tau \circ (\mu _B\ot A)=(A\ot \mu _B)\circ (\tau \ot B)\circ (B\ot
\tau).
\label{tw3}
\end{eqnarray}
If $\tau$ is a twisting map, then the map defined on $A\ot B$ by $\mu_\tau:=(\mu
_A\otimes \mu _B)\circ (A\otimes \tau \otimes B)$, is an associative
multiplication and $1\ot 1$ is the unit. The algebra having $A\ot B$ as
underlying vector space and $\mu_\tau$ as a product is denoted by
$A\ot_\tau B$ and is called the \dtext{twisted tensor product} of $A$ and $B$
(with respect to the twisting map $\tau$). For unital algebras, the existence
of a twisting map is equivalent to the existence of an algebra structure on
$A\ot B$ which is compatible with the canonical inclusions of $A$ and $B$, and
thus twisting maps are in one to one correspondence with the so-called
\dtext{factorization structures} introduced in \cite{majid}. From a categorical
point of view, a twisted tensor product is a particular situation of the theory
of \dtext{distributive laws} introduced by J. Beck in \cite{beck}.

Very little is known about the classification of the existing twisting maps
between two given algebras. Even in the simplest cases, this turns out to be a
very difficult problem to tackle. In \cite{cibils}, C. Cibils proposed a method
for describing all the twisting maps between $A$ and $B$, being $A=k^n$ the
algebra of functions over an $n$--points set, and $B=k^2$ the two-points
algebra. The resulting twisted tensor product algebras, which are dubbed
\dtext{noncommutative duplicates} can be realized up to some extent as a sort
of Ore extensions associated to the quotient algebra $k[x]/(x^2-x)$. For the
sake
of completeness, we sketch the procedure followed by Cibils for obtaining the
classification of the noncommutative duplicates.

\begin{proposition}
    The set of twisting maps between $A$ and $k^2$ (also called the set of
\dtext{2--interlacings of $A$}) is in one to one correspondence with the set
$Y_A$ of couples $(f,\delta)$ with $f\in\End{A}$ an algebra endomorphism and
$\delta:A\to A$ an idempotent $f$-twisted derivation such that
    \[  f=f^2+\delta f+ f\delta
    \]
\end{proposition}

Every algebra endomorphism $f$ of the algebra $A=k^n$ may be given
in terms of a set map $\varphi$, to which we can associate a
one-valued quiver with $n$ vertices. To this quiver, using the
derivation $\delta$ we may assign a \dtext{coloration} satisfying
certain conditions. Conversely, every one valued quiver which admits
a coloration satisfying those properties give rise to an algebra
endomorphism and a derivation as in the former proposition, and thus
to a twisting map, and so there is a one to one correspondence
between the set of 2-interlacings of $k^n$ and the set of coloured
one valued quivers with $n$ vertices.

Using this equivalence, Cibils gives a classification of all the
noncommutative duplicates of the algebras $k^n$, and computes their
Hochschild (co)homology using the techniques developed in
\cite{cibils98}. More concretely, the following results are used in
order to compute the Hochschild cohomology:

\begin{theorem}[\cite{cibils}]\label{thm_cibils}
Let $Q$ be a connected quiver which is not a crown, and let us denote by
$(KQ)_2$ the quotient of the path algebra $kQ$ by the two sided ideal
$(Q_{\geq 2})$ generated by the paths of length 2, then we have:
    \begin{enumerate}
    \item $\Dim_k HH^0 ((kQ)_2)= \#(Q_1 /\!\!/ Q_0) + 1$,
    \item $\Dim_k HH^1 ((kQ)_2)= \#(Q_1 /\!\!/ Q_1) -\#(Q_0) + 1$,
    \item $\Dim_k HH^n((kQ)_2) = \#(Q_n /\!\!/Q_1) - \#(Q_{n-1}/\!\!/Q_0)$ for
all $n\geq 2$.
    \end{enumerate}
    where for two sets of paths $X$ and $Y$, by $X /\!\!/ Y$ we denote the set
of \dtext{parallel paths}, that is, the set of couples $(x,y)\in X\times Y$
where $x$ and $y$ have the same source and target.
\end{theorem}

\begin{proposition}[\cite{cibils98}]\label{prop_cibils}
    Let $Q$ be a $c$--crown, with $c\geq 2$, then the center of $(kQ)_2$ is
one-dimensional. If the characteristic of $k$ is different from $2$,
for any $n$ which is an even multiple of $c$ we have
    \[
        \Dim_k HH^n \left((kQ)_2\right) = \Dim_{k}
HH^{n+1}\left((kQ)_2\right)=1.
    \]
    The cohomology vanishes in all other degrees.
\end{proposition}

These two results have an important consequence (see Corollary 3.2 of
\cite{cibils98}):

\begin{corollary}
    Let $Q$ be a connected quiver which is not a crown, then the graded
cohomology $HH^{\bullet}((kQ)_2)$ is finite dimensional if, and only if, $Q$
has no oriented cycles.
\end{corollary}


\section{The space of twisting maps}

Let us consider the algebras $A$ and $B$ both isomorphic to
$k[\mathbb{Z}_2]$, the ring algebra of the cyclic group
$\mathbb{Z}_2$, and let us fix $\langle 1_A,a \rangle$ basis of $A$
and $\langle 1_B,b \rangle$ basis of $B$, satisfying $a^2=1_A$ and
$b^2=1_B$. Then the set
    \[
        \langle 1_A \otimes 1_B, 1_A\otimes b,  a\otimes 1_B, a\otimes
        b \rangle
    \]
is a  basis of $A\ot B$, and
    \[
        \langle 1_B \otimes 1_A, b\otimes 1_A, 1_B\otimes a, b\otimes
        a \rangle
    \]
is a basis of $B\otimes A$.

The choice of these bases simplifies the required computations for
finding out all the twisting maps $\tau: B\otimes A \rightarrow A
\otimes B$, since the unitality conditions on $\tau$ forces us to
take
    \[
        \tau(1\ot 1)=1\ot 1,\ \tau(1\ot a)= a\ot 1,\ \tau(b\ot 1)=1\ot b,
    \]
so in order to give a twisting map between $A$ and $B$ it is enough
to give a value for $\tau(b\ot a)$ and check that it satisfies the
required compatibility conditions with respect to multiplications in
$A$ and $B$. In \cite{caenepeel}, an explicit approach to this
problem is performed, obtaining that any twisting map is one of the
following list:
\begin{enumerate}[$(a)$]
\item If $\Char(k)=2$, then:

  \begin{enumerate}[$(i)$]
   \item $\tau(b\otimes a)= \alpha (1_A \otimes 1_B)+ (a\otimes b)$, where
$\alpha \in k$.
   \item $\tau(b\otimes a)= \alpha (1_A \otimes 1_B)+ \alpha (1_A\otimes b)+
    \alpha (a\otimes 1_B)+ (\alpha+1) (a\otimes b)$, where $\alpha \in k$.
  \end{enumerate}
\item  If $\Char(k)\neq 2$, then:
  \begin{enumerate}[$(i)$]
  \item $\tau(b\otimes a)=(a\otimes b)$.
  \item $\tau(b\otimes a)=-(1_A\otimes 1_B)+\alpha (a\otimes b)$, where
$\alpha\in k$.
  \item $\tau(b\otimes a)= -(1_A \otimes 1_B)+(1_A \otimes b)+(a\otimes 1_B)$.
  \item $\tau(b\otimes a)= (1_A \otimes 1_B)- (1_A\otimes b)+ (a\otimes 1_B)$.
  \item $\tau(b\otimes a)= (1_A \otimes 1_B)+ \alpha (1_A\otimes b)- (a\otimes
1_B)$.
  \item $\tau(b\otimes a)= - (1_A \otimes 1_B)- (1_A\otimes b)- (a\otimes 1_B)$.
  \end{enumerate}
\end{enumerate}

The space of twisting maps over these particular algebras may also
be computed by means of certain coloured quivers, which are
associated to twisting maps following the procedure developed in
\cite{cibils}, as summarized at the end of Section
\ref{preliminaries}

In our situation, the algebra maps $f:k^2 \rightarrow k^2$ are all
given  as the lifting of the set maps $\varphi:\{a,b\} \rightarrow
\{a,b\}$ (see \cite{cibils} for this correspondence), thus obtaining
the four possible algebra maps given in generators by:
\begin{itemize}
\item $f_1(a)=a$ and $f_1(b)=b$.
\item $f_2(a)=b$ and $f_2(b)=a$.
\item $f_3(a)=a+b$ and $f_3(b)=0$.
\item $f_4(a)=0$ and $f_4(b)=a+b$.
\end{itemize}

Associated to these maps, we have the following quivers (where $Q_i$
stands for the quiver associated to the algebra map $f_i$):

\vspace{0.3cm}

$$\begin{array}{cc}
  Q_1:= \hspace{0.3cm}  \xymatrix{\circ \ar@(ur,ul)[] & \circ
\ar@(ur,ul)[]},\hspace{0.3cm}
   & Q_2:= \hspace{0.3cm} \xymatrix{\circ \ar@/^4pt/@<0.5 ex>[r]& \circ
\ar@/^4pt/@<0.5 ex>[l]}, \\
     & \\
     & \\
  Q_3:=\hspace{0.3cm} \xymatrix{\circ \ar[r] &  \circ
\ar@(ur,ul)[]},\hspace{0.3cm}

   & Q_4:=\hspace{0.3cm} \xymatrix{\circ\ar@(ur,ul)[] & \circ \ar[l]}
\end{array}$$

Now, the colorations attached to these quivers are given by:

\vspace{0.3cm}

\begin{enumerate}[($i'$)]
\item  \hspace{0.1cm} $\xymatrix{*+[o][F-]{\text{\scriptsize
$0$}} \ar@(ur,ul)[]& *+[o][F-]{\text{\scriptsize $0$}}
\ar@(ur,ul)[]}$

\vspace{0.2cm}

\item  $\xymatrix{*+[o][F-]{\text{\scriptsize
$\alpha$}} \ar@/^4pt/@<0.5 ex>[r]& *+[o][F-]{\text{\scriptsize
$\beta$}} \ar@/^4pt/@<0.5 ex>[l]}$ \hspace{0.05cm} where
$\beta=-1-\alpha$.

\vspace{0.2cm}

\item  \vspace{0.3cm} $\xymatrix{*+[o][F-]{\text{\tiny
$-1$}} \ar[r] &  *+[o][F-]{\text{\scriptsize $0$}} \ar@(ur,ul)[]}$,
 \hspace{0.2cm} $\xymatrix{*+[o][F-]{\text{\scriptsize
$0$}} \ar[r] &  *+[o][F-]{\text{\scriptsize $0$}} \ar@(ur,ul)[]}$,
 \hspace{0.2cm} $\xymatrix{*+[o][F-]{\text{\scriptsize
$0$}} \ar@(ur,ul)[] & *+[o][F-]{\text{\tiny $-1$}} \ar[l]}$ and
  \hspace{0.05cm} $\xymatrix{*+[o][F-]{\text{\scriptsize
$0$}} \ar@(ur,ul)[] & *+[o][F-]{\text{\scriptsize $0$}} \ar[l]}$
\end{enumerate}

Here we may observe that the twisting map ($i$) corresponds to the
coloured quiver ($i'$), the one--parameter family of maps ($ii$) is
associated to the quivers ($ii'$) when we vary the coloration, and
the twisting maps ($iii$), ($iv$), ($v$) and ($vi$) correspond to
the given colorations of ($iii'$).

\begin{remark} As a consequence of all this, the set of twisting maps gives rise
to a variety
consisting in five isolated points,  which correspond to the
twisting maps ($i$) and ($iii$)-($vi$), plus a $k$--line, associated
to the one-parameter family of maps described in ($ii$).
\end{remark}

\section{The isomorphism classes of the twisted algebras}

In the former section we described the set of all twisting maps between $k^2$
and $k^2$ but, as we mentioned earlier, different twisting maps could give rise
to isomorphic algebras. In this section we will describe the algebras
associated to the twisting maps that we obtained in the previous subsection,
describing the different isomorphism classes and giving a description of the
orbitspace in the corresponding variety of twisting maps.

In \cite{caenepeel}, a description of these algebras by means of generators and
relations is given, in particular mentioning that the algebras obtained from
the (non-invertible) twisting maps ($iii$)-($vi$) are all isomorphic to
    \[
        k\langle a,b\ |\ a^2=b^2=1,\ ba=a+b+1 \rangle.
    \]
    A different, but equivalent, description may be given following
\cite{cibils}, where it is shown that the algebras associated to the
four non-invertible twisting maps are all isomorphic to the path
algebra of the quiver
    \[  \widetilde{Q}:=
        \begin{array}{l}
        \xymatrix@R=10pt@C=10pt
        {   & \circ &  \\
            \circ \ar[rr] & & \circ
        }
        \end{array}
    \]
This means that four out of the five isolated points in our variety
provides the same point in the orbitspace. For the remaining
isolated point, which is the one corresponding to the flip map, i.e.
($i$), the corresponding algebra is just the usual tensor product:
    \[
        k\mathbb{Z}_2\ot k\mathbb{Z}_2\cong k\langle a,b\ |\ a^2=b^2=1,\
ba=ab\rangle.
    \]
Again, this algebra may be described as the path
algebra of the quiver
    \[
        \xymatrix@R=10pt@C=10pt
        {   \circ & \circ \\
            \circ & \circ
        }
    \]
This algebra is clearly non-isomorphic to the former one, since it
is commutative, and thus it gives a new point in the orbitspace.

Henceforth, the only remaining case is the one-parameter family of
twisting maps described in ($ii$). The family of algebras obtained
out of these twisting maps  is described in \cite{caenepeel} in
terms of generators and relations, obtaining the family
    \[
        A_q:=k\langle a,b\ |\ a^2=b^2=1,\ ab+ba=q\rangle,\quad \text{where $q\in
k$}.
    \]
    The authors of \cite{caenepeel} are not concerned by the number of different isomorphism
     classes of algebras which are obtained according to different
      values of the parameter. On the other hand, according to \cite[Theorem
4.4]{cibils}, all these algebras should be isomorphic to the
quotient of the path algebra of the so-called round-trip quiver
    \[  Q :=
        \xymatrix{
            \circ \ar@/^4pt/@<0.5ex>[r]& \circ \ar@/^4pt/@<0.5 ex>[l]
        }
    \]
    modulo the ideal generated by the set $Q_{\geq 2}$ of paths of length
greater
    than one. In other words, the obtained algebra would not depend on
the coloration. Unfortunately, the proof contains a slight mistake.
Within this one-parameter family of algebras we can find two
different kinds of algebras:
    \begin{itemize}
    \item If we take $q\neq \pm 2$, then the algebra map
        \[
            k\langle a,b\rangle \longrightarrow \mathcal{M}_2(k)
        \]
        defined by
        \[
            a\longmapsto \left(
                \begin{array}{cc}
                    1 & 0 \\
                    0 & -1
                \end{array} \right),\quad
            b\longmapsto \left(
                \begin{array}{cc}
                    \frac{q}{2} & \frac{2-q}{4} \\
                    \frac{2+q}{4} & -\frac{q}{2}
                \end{array}  \right)
        \]
provides an isomorphism of algebras between the algebra $A_q$ and
the $2\times 2$ matrix ring $\mathcal{M}_2(k)$.
    \item If $q\in \{2, -2\}$, the algebra map
    $f:A_{-2}\rightarrow A_2$ defined by
    \begin{gather*}
        f(1_A\otimes 1_B):=(1_A\otimes 1_B),\quad f(1_A\otimes
b):=(a\otimes 1_B), \\
        f(a\otimes 1_B):=(1_A\otimes b)-2 (a\otimes 1_B),\quad
f(a\otimes b):=-(a\otimes b)
    \end{gather*}
    is an isomorphism.

    Now, consider $R:=kQ/(Q_{\geq2})$ the quotient of the path algebra
of the round-trip quiver modulo the ideal  generated by $Q_{\geq
2}$. We may explicitly describe $R$ as the algebra having a basis
consisting in the four elements ${e, f, x, y}$ such that the
multiplication is given by the following table:
    \begin{center}
        \begin{tabular}{c|cccc}
        & $e$ & $f$ & $x$ & $y$ \\\hline
        $e$ & $e$ & 0 & 0 & $y$ \\
        $f$ & 0 & $f$ & $x$ & 0 \\
        $x$ & $x$ & 0 & 0 & 0  \\
        $y$ & 0  &  $y$ & 0  & 0
        \end{tabular}
    \end{center}

    Considering the algebra map $\phi: R \rightarrow A_{-2}$ defined
by:
    \begin{gather*}
        \phi(e):=1/2((1_A\otimes 1_B)-(a\otimes 1_B)),\\
        \phi(f):=1/2((1_A\otimes 1_B)+(a\otimes 1_B)), \\
        \phi(x):=1/4((1_A\otimes 1_B)+(1_A\otimes b)+(a\otimes 1_B)+(a\otimes
b)),\\
        \phi(y):=1/4((1_A\otimes 1_B)-(1_A\otimes b)-(a\otimes 1_B)+(a\otimes
b)),
    \end{gather*}
    we have $\phi$ is an algebra isomorphism between $A_{-2}$
and $R$, obtaining that both $A_2$ and $A_{-2}$ are isomorphic to
the algebra $R$.
    \end{itemize}
Finally, the line corresponding to our one-parameter family of
twisting maps corresponds to two points in the orbit space, one open
orbit corresponding to the matrix ring, plus one more point in the
closure of this orbit, corresponding to the quotient of the path
algebra of the round-trip quiver. From the point of view of
deformation theory, this means that the matrix algebra, realized as
a twisted tensor product, admits a deformation to the algebra
$kQ/(Q_{\geq 2})$.
    \begin{figure}[h]
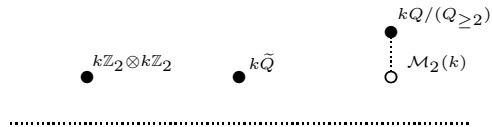

        \[
            \xy
                (10,6)*{\bullet}+(+6,2)*{\scriptscriptstyle k\mathbb{Z}_2\ot
                k\mathbb{Z}_2}; (30,6)*{\bullet}+(3,2)*{\scriptscriptstyle
                k\widetilde{Q}};
(50,6)*{\circ}+(6,2)*{\scriptscriptstyle \mathcal{M}_2(k)};
(50,12)*{\bullet}+(7,2)*{\scriptscriptstyle kQ/(Q_{\geq 2})};
                (0,0) \ar@{.} (60,0) 
                \POS (50,6)*{\circ} \ar@{.} (50,12);
            \endxy
        \]
    \caption{The orbit space}
    \end{figure}

Summarizing, we have proven the following result:

\begin{proposition}\label{classification}
 Let $k$ be a field with $\Char(k)\neq 2$. Let $A\cong B\cong k^2$, and let $\tau:B\ot A\to A\ot B$ be a twisting map,
then the twisted tensor product algebra $R:A\otimes_{\tau} B$ must
be isomorphic to one of the following algebras:
    \begin{enumerate}
        \item[(I)] $k^4$, or equivalently,
    the path algebra of the quiver
    
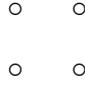
\[
        \xymatrix@R=10pt@C=10pt
        {   \circ & \circ \\
            \circ & \circ
        }
    \]
        \item[(IIa)] The algebra of matrices $\mathcal{M}_2(k)$.
        \item[(IIb)] The quotient $kQ/(Q_{\geq 2})$ of the path algebra $kQ$ of
the round-trip quiver
    \[  Q=
        
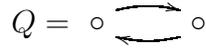
\xymatrix
        {   \circ \ar@/^4pt/@<0.5ex>[r]& \circ \ar@/^4pt/@<0.5 ex>[l]
        }
    \]
    modulo the ideal  generated by the set $Q_{\geq 2}$ of paths of length
greater than
one.
\item[(III)] The path algebra $k\widetilde{Q}$ of the quiver
   \[  \widetilde{Q}=
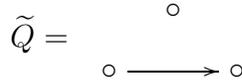
        \begin{array}{l}
        \xymatrix@R=10pt@C=10pt
        {   & \circ &  \\
            \circ \ar[rr] & & \circ
        }
        \end{array}
    \]
\end{enumerate}
\end{proposition}

\begin{remark} As we mentioned above, the classification given by Cibils for
noncommutative duplicates of set algebras in \cite{cibils}, is almost
complete, with the only exception being given by Theorem 4.4, dealing with
the connected components of the (coloured) quivers that are precisely the
round-trip quiver.
\end{remark}

Next we consider the formalism developed by Cibils for a two-fold
purpose, namely to highlight where the slight mistake in his proof
have been done and to obtain a characteristic free classification of
the isomorphism classes.  We have communicated to Cibils the
complete previous classification we have obtained, then he provided
us the precise localisation of the error in \cite{cibils}.

Following the same notation as Cibils does in \cite{cibils}, the
algebra structure of $A\otimes k[X]/(X^2-X)$ is determined by the
products
    \begin{equation}\label{product_cibils}
        Xa=\tau(X\otimes a)=\delta(a)+f(a)X
    \end{equation}
    for each $a\in A$, where $(\delta,f)$ is the
pair of the derivation and the endomorphism associated to the twisting map
$\tau$, see \cite[Proposition 2.10]{cibils}.

In our particular situation, that is, when we deal with the round-trip quiver,
the algebra endomorphism is given by
    \[ f(u)=v,\  f(v)=u
    \]
    whilst the derivation is given by
    \[ \delta(u)=a_v v-a_u u,\  \delta(v)=a_u u-a_v v,
    \]
    being $a_u$ and $a_v$ some parameters in $k$ and, $u$ and $v$ the
primitive orthogonal idempotent elements of $k^2=k\{u,v\}$ (cf.
\cite[Lemma 3.3]{cibils}). Applying formula \eqref{product_cibils}
to this particular situation we have:
    \begin{gather}
        Xu = -a_uu + a_uv + vX \\
        Xv = a_uu - a_uv + uX.
    \end{gather}
    Remember that in order to get a well-defined, associative structure, it is
necessary and sufficient to have $a_u+a_v+1=0$, as mentioned in
\cite[Theorem 3.14]{cibils}. Using this, the multiplication of the
resulting algebra may be summarized in the following table:
 \begin{center}
 \begin{tabular}{c|cccc}

     & $u$ & $uX$ & $v$ & $vX$ \\
    \hline
    $u$ & $u$ & $uX$ & 0 & 0 \\
    $uX$ & $-a_u u$ & $-a_u uX$ & $a_u u+uX$ & $-a_v uX$ \\
    $v$ & 0 & 0 & $v$ & $vX$  \\
    $vX$ & $a_v v+vX$  &  $-a_u vX$ & $-a_v v$  & $-a_vvX$ \\

\end{tabular} \end{center}

Now, observe that we have
    \begin{gather}
        (vXu)(uXv) = (vX)(uX)v = a_ua_vv, \\
        (uXv)(vXu) = (uX)(vX)u = a_ua_vu,
    \end{gather}
and this products are zero if, and only if, $a_ua_v=0$, a condition which
 is equivalent to have $a_u=0$ and $a_v=-1$, or $a_u=-1$ and
$a_v=0$. In this two cases we may carry on with the proof of
\cite[Theorem 4.4]{cibils}, obtaining the isomorphism with the quotient of the
path algebra of the round-trip quiver, as Cibils states (in our
classification these algebras correspond to $A_2$ and $A_{-2}$).

However, if the product $a_ua_v$ is non-zero, that is, if neither
$a_u$ nor $a_v$ are 0, then the map $\psi: kQ_f\rightarrow
k\{u,v\}\otimes k[X]/(X^2-X)$ considered in \cite[Theorem
4.4]{cibils} is no longer an algebra map. Still, for these cases it
is possible to consider the algebra isomorphism
    \begin{equation}
        f:k\{u,v\}\otimes k[X]/(X^2-X) \longrightarrow \mathcal{M}_2(k)
    \end{equation}
    given by
    \[
        \begin{array}{cc}
            u\longmapsto \left(
                \begin{array}{ll}
                    1 & 0 \\
                    0 & 0
                \end{array} \right),
            & v\longmapsto \left(
                \begin{array}{cc}
                    0 & 0 \\
                    0 & 1
                \end{array} \right), \\
            & \\
            uXv \longmapsto \left(
                \begin{array}{cc}
                    0 & a_ua_v \\
                    0 & 0
                \end{array} \right), &
            vXu \longmapsto \left(
                \begin{array}{cc}
                    0 & 0 \\
                    a_ua_v & 0
                \end{array} \right),
        \end{array}
    \]
in agreement with the result that we got in Proposition
\ref{classification}.

\begin{remark} It is worth noting that the fact that a matrix algebra may be
written in this way shows that the twisted tensor product of two
elementary algebras (as is the case for the algebras that we are
considering) is not in general a elementary algebra, even if we
require the twisting map to be bijective. Actually, the example that
we present shows that we can build a twisted tensor product of two
elementary algebras by means of an invertible twisting map, and
obtain an algebra which is not even basic!
\end{remark}

\section{Hochschild cohomology}

In this section we give a description of some facts related to the
Hochschild cohomology of the twisted tensor product algebras that we
have described above. Due to the similarity in the construction of
the twisted tensor product with the one performed for the usual
tensor product, it is reasonable to expect that Hochschild homology
groups should satisfy a sort of (maybe twisted) K\"{u}nneth formula
that would allow to compute the homology groups of the twisted
tensor product algebra out of the homology groups of the factors. A
step in this direction was given by J. A. Guccione and J. J.
Guccione in \cite{guccione}, where they build up a bicomplex which
should allow to compute the (co)homology for the twisted tensor
product when the twisting map is bijective, stating as a consequence
that the Hochschild dimension of a twisted tensor product is bounded
by the sum of the Hochschild dimensions of the factors. In
particular, this result would imply that any twisted tensor product
of two separable (i.e. having Hochschild dimension equal to 0)
algebras is again separable.

This result is false, and the counterexample we consider shows that
there is no hope to obtain a correct reformulation. We can build up
a twisted tensor product of two separable algebras (both of them
isomorphic to $k^2$) with respect to an invertible twisting map, and
such that the resulting algebra does not even have finite Hochschild
dimension. In order to do this, we give explicit descriptions, using
some methods developed by Cibils in \cite{cibils98} and
\cite{cibils}, of the Hochschild cohomology of all the algebras that
we classified in the former section.

\begin{proposition} Let $A\cong B\cong k^2$, and let $\tau:B\ot A\to A\ot B$ be
a twisting map, then the Hochschild cohomology of twisted tensor
product algebra $R:=A\otimes_{\tau} B$ is given by:

    \begin{enumerate}
        \item[(I)] If $R\cong k^4$, then $HH^0(R)=k^4$ and $HH^n(R)=0$ for any
$n\geq 0$.
        \item[(IIa)] If $R\cong \mathcal{M}_2(k)$, then  $HH^0(R)=k$ and
$HH^n(R)=0$ for any $n\geq 0$.
        \item[(IIb)] If $R\cong kQ/(Q_{\geq 2})$, then $HH^n(R)=k$
for all $n\geq 0$. In particular $R$ has infinite Hochschild
dimension.
        \item[(III)] If $R=k\widetilde{Q}$, then $HH^0(R)=k^3$, and
$HH^n(R)=0$ for all $n\geq 1$.
    \end{enumerate}
\end{proposition}

\begin{proof}
    The cases $(I)$ and $(IIa)$ are trivial, since both $k^4$ and $\mathcal{M}_2(k)$
are separable algebras (the latest because it is Morita equivalent
to the ground field $k$).

Case $(III)$ is a consequence of Theorem \ref{thm_cibils}.

Case $(IIb)$ is a direct consequence of Proposition
\ref{prop_cibils}. Since this is the situation that provides us the
aforementioned counterexample, for the sake of completeness, we
sketch Cibil's procedure applied to this particular example:

Recall (cf. \cite{cibils90a}, \cite{cibils90b}) that, if we have a
finite dimensional algebra $R$ admitting a decomposition $R=E\oplus
J$, being $E$ a maximal semisimple subalgebra of $R$ (which is
separable) and $J$ the Jacobson radical of $R$, then the Hochschild
cohomology of $R$ can be computed as the cohomology of the following
complex of cochains:
    \begin{equation}\label{complex}
        \xymatrix@C=11pt{ 0\ar[r] &  R^E \ar[r] & \Hom_{E-E}(J, R)\ar[r] & \ar@{.}[r]
        &
        \ar[r] & \Hom_{E-E}(J^{\ot^n_E}, R)\ar[r] & \ar@{.}[r] & }
    \end{equation}
%
    where $J^{\ot^n_E}$ is the tensor product over $E$ of $n$ copies of $J$.
    Whenever the Jacobson radical satisfies that $J^2=0$, the coboundary is given by
    \begin{gather*}
        (\delta r)(x):=rx-xr\quad\forall\, r\in R^E,\ x\in J,\\
        (\delta f)(x_1\ot\dotsb\ot x_{n+1}) = x_1 f(x_2\ot\dotsb\ot
        x_{n+1}) + (-1)^{n+1}f(x_1\ot\dotsb\ot x_n)x_{n+1}
    \end{gather*}
    for all $f\in \Hom_{E-E}(J^{\ot^n_E}, R)$.
    In our particular example, we have $kQ/(Q_{\geq 2})\cong
    kQ_0\oplus kQ_1$, being $E=kQ_0\cong k^2$ the (commutative)
    maximal semisimple subalgebra of $R$ and $kQ_1=J$ its Jacobson
    radical (whose square is 0). It is immediate to check that $J^{\ot^n_E}$ admits as a
    basis the set $Q_n$ of paths of length $n$. Now, using the additivity of the $\Hom$ functor, we have
    $\Hom_{E-E} (kQ_n, R)\cong \Hom_{E-E} (kQ_n, kQ_0)\op \Hom_{E-E} (kQ_n,
    kQ_1)$, and, as every simple subbimodule of $kQ_n$ corresponds to the
    bimodule generated by a path $\gamma$ of length $n$, which we can
    associate to the couple of vertices $(s(\gamma), t(\gamma))$ of
    starting and ending points of $\gamma$. Applying Schur's lemma,
    we have $\Hom_{E-E} (k\gamma, k\gamma')=0$ unless $\gamma$ and
    $\gamma'$ have the same starting and ending points, that is,
    unless $\gamma$ and $\gamma'$ are \dtext{parallel paths}. Using
    this, we find a linear isomorphism $\Hom_{E-E} (kQ_n, kQ_0)\simeq k(Q_n/\!\!/
    Q_0)$. Similarly, we have a linear isomorphism $\Hom_{E-E} (kQ_n, kQ_1)\simeq k(Q_n/\!\!/
    Q_1)$. Through these identifications, the coboundary $\delta$ is
    translated into the coboundary $\left( \begin{smallmatrix} 0 & 0 \\ D & 0
    \end{smallmatrix}\right)$, where the map $D:k(Q_n/\!\!/Q_0)\to~k(Q_{n+1}/\!\!/Q_1)$ is given by
    \[
        D(\gamma, e):=\sum_{a\in Q_1 e} (a\gamma, \gamma) +
        (-1)^{n+1}\sum_{a\in eQ_1} (\gamma a, a).
    \]
    By construction, we obtain a complex isomorphism between
    \eqref{complex} and the complex
    \begin{equation}\label{complex2}
        0\to kQ_0\lto k(Q_1/\!\!/Q_0)\oplus k(Q_1/\!\!/Q_1) \lto
        \dotsb \lto k(Q_n/\!\!/Q_0)\oplus k(Q_n/\!\!/Q_1) \lto
        \dotsb.
    \end{equation}
    Since our quiver has no loops, whenever $n$ is odd we have
    \[
        k(Q_n/\!\!/Q_0)=k(Q_{n+1}/\!\!/Q_1)=\{0\},
    \]
    whilst for $n$ even we get
    \[
        k(Q_n/\!\!/Q_0)\cong k(Q_{n+1}/\!\!/Q_1)\cong k\mathbb{Z}_2,
    \]
    as every path is uniquely determined by its starting (and
    ending) point, where the identification consists on sending a
    path to $1$ if it starts at the vertex $e$ or to $t$ if it
    starts at the vertex $f$, and we are considering $k\mathbb{Z}_2=k\{1,t|\ t^2=1\}$.
    Via this identification, for even $n$, the map $D$ transforms into the map
    $D':k\mathbb{Z}_2\to k\mathbb{Z}_2$ defined by
    \[ D'(1)=1-t,\ D'(t)=t-1.
    \]
    This map obviously has one dimensional kernel, generated by the
    element $1+t$, and one dimensional image. Summing everything up,
    we may rewrite the complex \eqref{complex2} as
    \[
        \xymatrix@C=15pt{
        0\ar[r]& k^2 \ar[rr]^{D'}&& k^2 \ar[rr]^{0} &&
        k^2\ar[rr]^{D'} &&
        \ar@{.}[r] & },
    \]
    and thus, for $n$ odd we have
    \[ \Dim_k HH^n\left((kQ)_2\right)=\Dim_k \left(\frac{\ker
    0}{\im D'}\right)=\Dim_k k^2 - \Dim_k (\im D')=1,
    \]
    whilst, for $n$ even we get
    \[ \Dim_k HH^n\left((kQ)_2\right)=\Dim_k \left(\frac{\ker
    D'}{\im 0}\right)=\Dim_k (\ker D') - \Dim_k (0)=1,
    \]
    as we wanted to prove.
\end{proof}

As we announced, the algebra of type ($IIb$) provides us an example
of a twisted tensor product of two separable algebras, with respect
to a bijective twisting map, which does not have finite Hochschild
dimension. This example contradicts \cite[Theorem 1.6, Theorem 1.7,
Corollary 1.8]{guccione}. It is worth noting that in order to
disprove Gucciones' results, it is not necessary to give an explicit
description of the Hochschild cohomology, being enough to show that
the twisted tensor product algebra is not separable. An immediate
prove of this fact follows from the realization of this algebra as
the quotient $R = kQ/(Q_{\geq 2})$, as we can immediately check that
the elements of $R$ corresponding to (the equivalence classes of)
the arrows of $Q$ provide nonzero elements of the Jacobson radical
of $R$ (actually, the Jacobson radical is precisely the ideal
generated by these two elements).


\section*{Acknowledgements}
 The authors would like to thank Dragos Stefan for helping us to
understand the subtleties beyond Hochschild cohomology, and Calude
Cibils, for his willingness to check our work and his numerous
remarks and corrections to the preliminary version of this paper.

\end{document}